\documentclass[11pt,a4paper]{article}
\usepackage[a4paper, margin=2.5cm]{geometry}

\usepackage{amsmath,amsthm}

\usepackage{hyperref}
\usepackage{dsfont}
\usepackage{amsmath}
\usepackage{enumitem}
\usepackage{amssymb}
\usepackage{amsfonts}
\usepackage{graphicx}
\usepackage{subfigure}
\usepackage{color}

\newcommand{\red}[1]{\textcolor{black}{#1}}

\newcommand{\argmax}{\operatornamewithlimits{argmax}}
\newcommand{\argmin}{\operatornamewithlimits{argmin}}

\theoremstyle{plain}
\newtheorem{theorem}{Theorem}[section]
\theoremstyle{plain}
\newtheorem{proposition}{Proposition}[section]
\theoremstyle{definition}
\newtheorem{definition}{Definition}[section]
\theoremstyle{corollary}
\newtheorem{corollary}{Corollary}[section]

\title{Revisiting maximum-a-posteriori estimation in log-concave models}

\author{Marcelo Pereyra \footnote{School of Mathematical and Computer Sciences, Heriot-Watt University (m.pereyra@hw.ac.uk)}}

\begin{document}

\maketitle

\begin{abstract}
Maximum-a-posteriori (MAP) estimation is the main Bayesian estimation methodology in imaging sciences, where high dimensionality is often addressed by using Bayesian models that are log-concave and whose posterior mode can be computed efficiently by convex optimisation. However, despite its success and wide adoption, MAP estimation is not theoretically well understood yet. In particular, the prevalent view in the community is that MAP estimation is not proper Bayesian estimation in the sense of Bayesian decision theory because it does not minimise a meaningful expected loss function (unlike the minimum mean squared error (MMSE) estimator that minimises the mean squared loss). This paper addresses this theoretical gap by presenting a general decision-theoretic derivation of MAP estimation in Bayesian models that are log-concave. A main novelty is that our analysis is based on differential geometry, and proceeds as follows. First, we use the underlying convex geometry of the Bayesian model to induce a Riemannian geometry on the parameter space. We then use differential geometry to identify the so-called natural or canonical loss function to perform Bayesian point estimation in that Riemannian manifold. For log-concave models, this canonical loss coincides with the Bregman divergence associated with the negative log posterior density. Following on from this, we show that the MAP estimator is the only Bayesian estimator that minimises the expected canonical loss, and that the posterior mean or MMSE estimator minimises the dual canonical loss. We then study the question of MAP and MSSE estimation performance in high dimensions. Precisely, we establish a universal bound on the expected canonical error as a function of image dimension, providing new insights the good empirical performance observed in convex problems. Together, these results provide a new understanding of MAP and MMSE estimation in log-concave settings, and of the multiple beneficial roles that convex geometry plays in imaging problems. {Finally, we illustrate this new theory by analysing the regularisation-by-denoising Bayesian models, a class of state-of-the-art imaging models where priors are defined implicitly through image denoising algorithms, and an image denoising model with a wavelet shrinkage prior.}
\end{abstract}

\section{Introduction}\label{sec:BayesianPointEstimation}
We consider the estimation of an unknown image $x \in \mathbb{R}^n$ from some data $y$, related to $x$ by a statistical model with likelihood $p(y|x)$. Adopting a Bayesian statistical approach, we postulate a prior distribution $p(x)$ modelling the prior knowledge available, and base our inferences on the posterior distribution \cite{somersalo:2005}
$$
p(x|y) = \frac{p(y|x)p(x)}{\int_{\mathbb{R}^n} p(y|x)p(x) \textrm{d}x},
$$
which models our knowledge about $x$ after observing $y$. In this paper we focus on the case where $p(x|y)$ belongs to the class of log-concave distribution, i.e.,
\begin{equation}\label{posterior}
p(x|y) = \frac{\exp\{-\phi(x)\}}{\int_{\mathbb{R}^n} \exp\{-\phi(s)\} \textnormal{d}s},
\end{equation}
for some proper convex function $\phi: \mathbb{R}^n \rightarrow (-\infty,\infty]$ \cite{CombettesBook}, and where we note that for notation convenience we do not write explicitly the dependence of $\phi$ on $y$.

Log-concave models \eqref{posterior} are ubiquitous in modern imaging sciences. For instance, many imaging methods to solve linear imaging inverse problems involving additive Gaussian noise use models of the form $\phi(x) = \|y-Ax\|^2 /2\sigma^2 + \psi(Bx) + \boldsymbol{1}_\mathcal{S}(x)$ for some linear operators $A$ and $B$, convex regulariser $\psi$, and convex set $\mathcal{S}$ (see \cite{Figueiredo2011,Babacan:2010,Babacan:2011} for examples related to image deblurring, inpainting, compressive sensing, super-resolution, and tomographic reconstruction, with total-variation and wavelet priors). {Similar log-concave Bayesian models can be considered for problems involving other observation noise models, such as Poisson noise \cite{PIDAL:2010}, and for phase retrieval problems \cite{PUMA:2007}}. Log-concave models \eqref{posterior} are also used extensively in other areas of data science such as machine learning \cite{theodoris}.

Because drawing conclusions directly from $p(x|y)$ is difficult, Bayesian imaging methods deliver summaries of $p(x|y)$ - namely Bayes point estimators - that summarise $p(x|y)$ optimally in a single value $\hat{x}$. This estimated value is optimal the following decision-theoretic sense \cite{cprbayes}:\vspace{-0.2cm}\\

\begin{definition}\label{defBayesEst}
Let $L : \mathbb{R}^n \times \mathbb{R}^n \rightarrow \mathbb{R}_0^+$ be a loss function that quantifies the difference between two points in $\mathbb{R}^n$. A Bayes estimator associated with $L$ is any estimator that minimises the posterior expected loss, i.e.,
\begin{equation*}
\begin{split}
\hat{x}_L &= \argmin_{u \in \mathbb{R}^n} \textnormal{E}_{x|y}[L(u,x)].
\end{split}
\end{equation*}
\end{definition}
{Recall that the posterior expectation $\textnormal{E}_{x|y}[L(u,x)] \triangleq \int_{\mathbb{R}^n} L(u,x) p(x|y) \textnormal{d}x$}. Sensible loss functions $L$ usually verify the following three desiderata \cite{cprbayes}:\vspace{0.2cm}

\begin{itemize}[leftmargin=1cm]
\item $L(u,x) \geq 0,\,\forall \,u, x \in \mathbb{R}^n$\,,
\item $L(u,x) = 0\,\, \iff \,\,u = x$\,,  
\item $L$ strictly convex w.r.t. its first argument (to guarantee estimator uniqueness).
\end{itemize}
\vspace{0.2cm}
\noindent {Estimator uniqueness is important because it implies admissibility (i.e., Bayesian estimator $\hat{x}_L$ is not dominated by any other estimator) \cite{cprbayes}}. Observe that $L$ is not necessarily symmetric; that is, $L(u,x) \neq L(x,u)$. We do not enforce symmetry because the arguments of $L$ have clearly different roles in the decision problem. 

In a purely theoretical Bayesian exercise, $L$ should be chosen based on specific aspects of the problem and application considered. This is particularly important in imaging problems that are ill-posed or ill-conditioned, as the choice of $L$ can significantly impact estimation results. However, specifying a bespoke loss function for high dimensional problems is not easy. Consequently, most methods in the imaging literature use default losses and estimators.

In particular, Bayesian imaging methods have traditionally used the minimum mean squared error (MMSE) estimator, {given by the posterior mean $\hat{x}_{\text{MMSE}} =  \int_{\mathbb{R}^n} p(x|y)\,x \textnormal{d}x$}. This estimator is widely regarded as a gold standard in the field, in part because of its good empirical performance and favourable theoretical properties, and also perhaps in part  because of cultural heritage. {From Bayesian decision theory, MMSE estimation is optimal with respect to the entire class of Euclidean or Mahalanobis squared distances, given by quadratic loss functions of the form $L(u,x) = (u-x)^\top Q (u-x)$ with $Q \in \mathcal{S}^n_{++}$ (i.e., the set of $n\times n$ positive definite matrices), which includes the popular mean square loss $L(u,x) = \|u-x\|_2^2$ when $Q = I_n$ \cite{cprbayes}}. This gives $\hat{x}_{\text{MMSE}}$ a straightforward geometric interpretation. Moreover, MSE estimation is optimal w.r.t. a more general class of functions \cite{Banerjee:2005} that provides a second order approximation to any strongly convex loss function; hence $\hat{x}_{\text{MMSE}}$ can act as a universal proxy for other Bayesian estimators in this sense. 


Unfortunately, calculating $\hat{x}_{\text{MMSE}}$ in high dimensional models can be very difficult because it requires solving integrals that are often too computationally expensive. This has stimulated much research on the topic, from fast Monte Carlo simulation methods to efficient approximations with deterministic algorithms \cite{Pereyra2016,durmus2018}. But with ever increasingly large problems and datasets, many imaging methods have focused on alternatives to MMSE estimation.

In particular, modern methods rely strongly on maximum-a-posteriori (MAP) estimation
\begin{align*}
\hat{x}_{\text{MAP}} &= \argmax_{x \in \mathbb{R}^n} p(x|y),\\
&= \argmin_{x \in \mathbb{R}^n} \phi(x),
\end{align*}
whose calculation is a convex problem that can often be solved very efficiently, even in very high dimensions (e.g., $n > 10^6$), by using proximal convex optimisation techniques \cite{PockActaNumerica,boydbook,Combettes2011,Green2015}. Modern non-statistical imaging methods also predominately solve problems by convex optimisation, and their solutions are often equivalent to performing MAP estimation for some implicit Bayesian model. The precise sense in which these solutions are equivalent to MAP estimators is an interesting discussion topic that is beyond the scope of this paper.

Following a decade of intensive activity, there is now abundant evidence that MAP estimation delivers accurate results for a wide range of imaging problems. However, from a theoretical viewpoint MAP estimation is not well understood. Currently the predominant view is that MAP estimation is not formal Bayesian estimation in the decision-theoretic sense postulated by Definition \ref{defBayesEst} because it does not generally minimise a known expected loss. The prevailing interpretation is that MAP estimation is in fact an approximation arising from the degenerate loss $L_\epsilon(u,x) = \boldsymbol{1}_{\|x-u\|<\epsilon}$ with $\epsilon \rightarrow 0$ \cite{cprbayes} (this derivation holds for all log-concave models, but is not generally true \cite{Bassett2016}). However, this asymptotic derivation does not lead to a proper Bayesian estimator. More importantly, the resulting loss is very difficult to motivate for inference problems in continuous domains such as $\mathbb{R}^n$, and does not help explain the good empirical performance reported in the literature. 

Furthermore, most other theoretical results for MAP estimation only hold for very specific models, or have been derived by adopting analyses that are extrinsic to the Bayesian decision theory framework (e.g. by analysing MAP estimation as constrained or regularised least-squares regression, see for example \cite{Candes2006,Candes2009}). As a trivial example of results that only hold for specific models, when $p(x|y)$ is symmetric we have $\hat{x}_{\text{MAP}} = \hat{x}_{\text{MMSE}}$, and thus MAP estimation inherits the favourable properties of MMSE estimation. This result has been partially extended to some denoising models of the form $p(x|y) \propto \exp\{\|y-x\|^2/2\sigma^2 + \lambda h(x)\}$ in \cite{Gribonval:2011}, where it is shown that MAP estimation coincides with MMSE estimation with a different model $\tilde{p}(x|y) \propto \exp\{- \|y-x\|^2/2\sigma^2 - \tilde{\lambda} \tilde{h}(x)\}$ involving a different prior distribution. It follows that for these models MAP estimation is decision-theoretic Bayesian estimation w.r.t. the weighted loss $L(u,x) = \|u-x\| \exp\{\tilde{\lambda}\tilde{h}(x)-\lambda h(x)\}$. This is of course a post-hoc loss, but the result is interesting in that it highlights that a single estimator may have a plurality of origins. More importantly, it raises the question if MAP estimation is merely a computational proxy for MMSE estimation, which unlike $\hat{x}_{\text{MAP}}$ has an appealing theoretical underpinning. {This question was recently answered in Burger \& Lucka \cite{Lucka:2014}: MAP estimation is proper decision-theoretic Bayesian estimation for all models of the form $p(x|y) \propto \exp\{-\|y-Ax\|^2_{\Sigma^{-1}}/2 - \lambda h(x)\}$, with known linear operator $A$ and noise covariance $\Sigma$, and where $h$ is convex and Lipchitz continuous. More precisely, that paper shows that for this class of models MAP estimation is optimal w.r.t. the loss $L(u,x) = \|A(u-x)\|^2_{\Sigma^{-1}} + 2\lambda D_h(u,x)\}$, where $D_h(x) = h(u) - h(x) - \nabla h(x)^\top(u-x)$ is the $h$-Bregman divergence \cite{CombettesBook}. The paper also shows that $\hat{x}_{\text{MAP}}$ outperforms $\hat{x}_{\text{MMSE}}$ w.r.t., the expected Bregman error $\textnormal{E}_{x|y} \{D_h(u,x)\}$, an error measure that grasps some distinctive features of $x$ (e.g., sparsity, regularity, smoothness, etc.). It may appear that the loss identified in \cite{Lucka:2014} is rather artificial and difficult to analyse and motivate; however, the new results presented in Section \ref{sec:CanonicalBayesian} show that it is a specific instance of a more general loss that stems directly from the consideration of the geometry of the Bayesian model.}

{It is also worth mentioning that several recent works have studied MAP estimation in infinite-dimensional settings, which is important for our understanding of how the technique behaves in increasingly large problems. An important advance in this area is the connection of the topological description of the MAP estimate to a variational problem, developed in \cite{Dashti} for non-linear inverse problems in a Gaussian framework, and subsequently extended to non-Gaussian settings in \cite{Helin}. Agapiou et al. made another key contributions in this area by studying infinite-dimensional MAP estimation with Besov priors, which are very relevant to imaging sciences because they promote sparsity and preservation of edges\cite{Agapiou}. We also mention the recent work \cite{Sullivan} that further improves our understanding of modes in infinite dimensions.}

In order to better understand MAP estimation, in this paper we first revisit the choice of the loss function for Bayesian point estimation in the context of models that are log-concave, where MAP is a convex problem (we limit our analysis to finite-dimensional problems). A main novelty of our analysis is that, instead of specifying the loss directly, we use differential geometry to derive the loss from the geometry of the model. Precisely, we show that under some regularity assumptions, the log-concavity of $p(x|y)$ induces a specific Riemannian differential geometry on the parameter space, and that taking into account this space geometry naturally leads to an intrinsic or {canonical} loss function to perform Bayesian point estimation in that space. Following on from this, we establish that the canonical loss for the parameter space is the Bregman divergence associated with $\phi(x) = -\log p(x|y)$, and that the Bayesian estimator w.r.t. this loss is the MAP estimator. We then show that the MMSE estimator is the Bayesian estimator associated with the dual canonical loss, and propose universal estimation performance guarantees for MAP and MMSE estimation in log-concave models. {We conclude by illustrating our theory with an application to linear inverse problems with sparsity-promoting wavelet priors, and an analysis of the regularisation-by-denoising models proposed recently in \cite{RED:2017}.}

The remainder of the paper is organised as follows: Section \ref{sec:Riemannian} introduces the elements of differential geometry that are essential to our analysis. In Section \ref{sec:CanonicalBayesian} we present our main theoretical results: a decision-theoretic and differential-geometric derivation of MAP and MMSE estimation, with universal bounds on the estimation error involved. Section \ref{relax} discusses the impact of relaxing the regularity assumptions adopted in Section \ref{sec:CanonicalBayesian}. Conclusions are finally reported in section \ref{sec:Conclusion}. Proofs are presented in the appendix.

\section{Riemannian geometry and the canonical divergence function} \label{sec:Riemannian}
In this section we recall some elements of differential geometry that are necessary for our analysis. For a detailed introduction to this topic we refer the reader to \cite{amaribook}.

{
An $n$-dimensional Riemannian manifold $(\mathbb{R}^n,g)$, with metric $g : \mathbb{R}^n \rightarrow \mathcal{S}^n_{++}$ and global coordinate system $x$, is a vector space that behaves locally as an Euclidean space\footnote{Recall that $\mathcal{S}^n_{++}$ is the set of $n\times n$ positive definite matrices.}. More precisely, at any point $x \in \mathbb{R}^n$, we have a tangent space $\mathcal{T}_x \mathbb{R}^n$ with inner product $\langle u,x\rangle = u^\top g(x) x$ and norm $\|x\| =\sqrt{x^\top g(x) x}$. This tangent space describes how the manifold $(\mathbb{R}^n,g)$ behaves locally at $x$. The geometry is local and may vary smoothly from $\mathcal{T}_x \mathbb{R}^n$ to neighbouring tangent spaces (i.e., the inner product and norm used are local properties and vary spatially). The variations are encoded in the affine connection $\Gamma$, with coefficients given by $\Gamma_{i j,\,k} (x) = \partial_k g_{i,j}(x)$ describing the spatial evolution of the metric $g$.} 

A crucial property of $(\mathbb{R}^n,g)$ is that, similarly to Euclidean spaces, manifolds supports divergence functions:\vspace{-0.2cm}\\

\begin{definition}\label{divergences}
A function $D : \mathbb{R}^n \times \mathbb{R}^n \rightarrow \mathbb{R}$ is a divergence function on $\mathbb{R}^n$ if the following conditions hold for any $u,x \in \mathbb{R}^n$:
\begin{itemize}[leftmargin=1cm]
\item $D(u,x) \geq 0, \,\forall \,u, x \in \mathbb{R}^n$,
\item $D(u,x) = 0\,\, \iff \,\,x = u$,
\item $D(u,x)$ is strongly convex w.r.t. $u$, and $\mathcal{C}^2$ w.r.t $u$ and $x$.
\end{itemize}
\end{definition}

\vspace{0.2cm}Observe that the class of divergence functions is equivalent to that of loss functions for Bayesian point estimation specified in Section \ref{sec:BayesianPointEstimation}, with some mild additional regularity conditions. This suggests that divergence functions are sensible losses to define estimators. Moreover, divergence functions also provide a link to the differential geometry of the space, which allows relating space geometry and Bayesian decision theory. This relationship has been used previously to analyse Bayesian decision problems from a Riemannian geometric viewpoint, leading to the so-called decision geometry framework \cite{Dawid2007}. Here we adopt an opposite perspective: we start by considering a Riemannian manifold $(\mathbb{R}^n,g)$ and then use the relationship to identify the divergence functions that arise naturally in that space. In particular, we focus on the so-called {canonical} divergence on $(\mathbb{R}^n, g)$, which generalises the Euclidean squared distance to this kind of manifold \cite{amari:2015}.\vspace{0.2cm}

\begin{definition}[Canonical divergence \cite{amari:2015}]
For any two points $u,x \in \mathbb{R}^n$, the $(\mathbb{R}^n,g)$-canonical divergence is given by
\begin{equation}\label{canonicalDiv}
D(u,x) = \int_0^1 t\dot{\gamma_t}^\top g(\gamma_t) \dot{\gamma}_t \textnormal{d}t
\end{equation}
where $\gamma_t$ is the $\Gamma$-geodesic from $u$ to $x$ and $\dot{\gamma}_t = {\textnormal{d}}/{\textnormal{d}t} \,\gamma_t$.
\end{definition}
{The reason why $D$ is the $(\mathbb{R}^n,g)$-canonical divergence is that it fully specifies the geometry of $(\mathbb{R}^n,g)$; i.e., there is a one-to-one relationship between $D$ and the metric $g$.}

{Observe that $D$ has connections to the length of the $\Gamma$-geodesic between $u$ and $x$. Precisely, by noting that the squared length of a curve $\zeta_t : [0,1] \rightarrow \mathbb{R}^n$ on the manifold $(\mathbb{R}^n,g)$ is given by $ \int_0^1 \dot{\zeta_t}^\top g(\zeta_t) \dot{\zeta}_t \textnormal{d}t$, we observe that $D(u,x)$ is essentially the squared length of the $\Gamma$-geodesic $\gamma_t$ weighted linearly along the path from $u$ to $x$. This weighting in \eqref{canonicalDiv} is important because it guarantees that $D(u,x)$ is convex in $u$, a necessary condition to define a divergence function (the weighting also leads to other important properties such as linearity w.r.t. $g$, see section \ref{sec:CanonicalBayesian}).} Note that the weighting also introduces an asymmetry, i.e., generally $D(u,x) \neq D(x,u)$, which will have deep implications for Bayesian estimation.

Finally, it is easy to check that \eqref{canonicalDiv} reduces to the Euclidean squared distance $D(u,x) = \tfrac{1}{2} (u-x)^\top g (u-x)$ when $(\mathbb{R}^n, g)$ is the Euclidean space with product $\langle u,x\rangle = u^\top g x$\footnote{In the Euclidean case we have that case $g$ is constant, the $\Gamma$-geodesic is $\gamma_t = u + t(x-u)$, so $D(u,x) = \int_0^1 t (u-x)^\top g (u-x) \textnormal{d}t = \int_0^1 t \textnormal{d}t (x-u)^\top g (x-u) = \tfrac{1}{2} (x-u)^\top g (x-u)$.}. More generally, $D$ is always consistent with the local Euclidean geometry of the manifold $(\mathbb{R}^n, g)$. That is, for any point $x + dx$ in the neighbourhood of $x$ we have $D(x+dx,x) = \|dx\|^2/2 + {o}(\|dx\|^2)$, where $\|\cdot\|$ is the Euclidean norm of the tangent space $\mathcal{T}_x \mathbb{R}^n$ (a higher order approximation of $D(x+dx,x)$ is also possible by using the connection $\Gamma$ \cite{amaribook}). {And because $D$ is the canonical divergence, if we use the decision geometry framework \cite{Dawid2007} to derive the Riemannian geometry induced by $D$ on $\mathbb{R}^n$ we obtain
$$
g^{(D)}_{i,j}(x) \triangleq  \partial_i \partial_j D(x,x) = g_{i,j}(x),
$$
$$
\Gamma_{i j,\,k}^{(D)}(x) \triangleq \partial_i \partial_j \partial^\prime_k D(x,x) = \Gamma_{i j,\,k}(x),
$$ 
(here $\partial$ and $\partial^\prime$ denote differentiation w.r.t. the first and second components of $D$ respectively), indicating that $D$ fully specifies the geometry of $(\mathbb{R}^n,g)$, and vice-versa.}

\section{A differential-geometric derivation of MAP and MMSE estimation}\label{sec:CanonicalBayesian}
\subsection{Canonical Bayesian estimation: from differential geometry to decision theory}

In this section we use differential geometry to relate the geometry of $p(x|y)$ to the loss functions used for Bayesian estimation of $x$. Precisely, we exploit the log-concavity of $p(x|y)$ to induce a Riemannian geometry on the solutions space. This in turn defines a canonical loss for that space and two Bayesian estimators: a primal estimator related to $D(u,x)$ and a dual estimator related to the dual divergence $D^*_\phi(u,x) = D_\phi(x,u)$. We focus on the case where $p(x|y)$ is smooth and strongly log-concave, and later analyse the effect of relaxing these assumptions.\vspace{-0.2cm}\\

\begin{theorem}[Canonical Bayesian estimators]\label{Theo2}
Suppose that $\phi(x) = -\log p(x|y)$ is strongly convex, continuous, and $\mathcal{C}^3$ on $\mathbb{R}^n$. Let $(\mathbb{R}^n, g)$ denote the Riemannian manifold induced by $\phi$, with metric coefficients $g_{i,j}(x) = \partial_i \partial_j \phi(x)$. Then, the canonical divergence on $(\mathbb{R}^n, g)$ is the $\phi$-Bregman divergence, i.e.,
$$
D_\phi (u,x) = \phi(u) - \phi(x) - \nabla \phi(x)(u - x).
$$
In addition, the Bayesian estimator associated with $D_\phi(u,x)$ is unique and is given by the maximum-a-posteriori estimator,
\begin{align*}
\hat{x}_{D_\phi} &\triangleq \argmin_{u \in \mathbb{R}^n} \textnormal{E}_{x|y}[D_\phi(u,x)]\,,\\
							&=  \argmin_{x \in \mathbb{R}^n} \phi(x)\,,\\
						  &=  \hat{x}_{\text{MAP}}\,.
\end{align*}
The Bayesian estimator associated with the dual canonical divergence $D^*_\phi (u,x) = D_\phi(x,u)$ is also unique and is given by the minimum mean squared error estimator
\begin{align*}
\hat{x}_{D^*_\phi} &\triangleq \argmin_{u \in \mathbb{R}^n} \textnormal{E}_{x|y}[D^*_\phi(u,x)]\,,\\
							&=  \int_{\mathbb{R}^n} x p(x|y) \textnormal{d}x\,,\\
						     &=  \hat{x}_{\text{MMSE}}\,.
\end{align*}
\end{theorem}
\noindent The proof is reported in the appendix.

Theorem \ref{Theo2} provides several interesting new insights into MAP and MMSE estimation in log-concave models. First and foremost, MAP estimation stems from Bayesian decision theory, and hence it stands on the same theoretical footing as the core Bayesian methodologies such as MMSE estimation (albeit w.r.t. a different class of loss functions). The MAP loss, $D_\phi (u,x)$, is a generalisation of the Euclidean squared distance that arises naturally from the consideration of the geometry of $p(x|y)$. {Consequently, the conventional definition of the MAP estimator as the maximiser $\hat{x}_{\text{MAP}} = \argmax_{x \in \mathbb{R}^n} p(x|y)$ is mainly algorithmic for these models, useful to highlight that these estimators take the form of a convex optimisation problem. (Of course, this is a key computational advantage over other Bayesian point estimators because it allows computing $\hat{x}_{\text{MAP}}$ by using using modern proximal convex optimisation algorithms that scale very efficiently to high-dimensions - see e.g., \cite{PockActaNumerica} for details)}. Moreover, Theorem \ref{Theo2} also reveals a surprising form of duality between MAP and MMSE estimation, with the two estimators intimately related to each other by the (asymmetry of the) canonical divergence function that $p(x|y)$ induced on the solutions space. Note that Gaussian models are particular because $(\mathbb{R}^n, g)$ is Euclidean in that case, which is a self-dual space; consequently $D_\phi(u,x) = D_\phi(x,u) = \tfrac{1}{2}\|u-x\|_{\Sigma^{-1}}^2$ and the primal and dual canonical estimators coincide as a result. Finally, Theorem \ref{Theo2} also shows that, under log-concavity and smoothness assumptions, the posterior mode is a global property of $p(x|y)$, similarly to the posterior mean.

The way in which the Bregman divergence $D_\phi (u,x)$ measures the similarity between $u$ and $x$ is directly related to the geometry of $p(x|y)$. Precisely, because $\phi(x) = - \log p(x|y)$ is strongly convex, then $\phi(u) > \phi(x) - \nabla \phi(x)(u - x)$ for any $u \neq x$. The divergence $D_\phi (u,x)$ essentially quantifies this gap, which as mentioned previously, is directly related to the length of the affine geodesic from $u$ to $x$ (and hence not only to the relative position of $u$ and $x$ but also to the space geometry induced by $p(x|y)$). Moreover, this geometry can depend on the value of $y$, however for the important class of models of the $p(x|y) \propto \exp\{-\|y-Ax\|^2_{\Sigma^{-1}}/2 - \lambda h(x)\}$ the geometry is completely specified by $\Sigma$ and $\lambda h$ independently of $y$. Furthermore, observe that because $D_\phi$ is linear w.r.t. $\phi$, then for any decomposition $\phi = \alpha \phi_1 + \beta \phi_2$ based on two convex functions $\phi_1$ and $\phi_2$ and $\alpha, \beta \in \mathbb{R}$, we obtain $D_\phi = \alpha D_{\phi_1} + \beta D_{\phi_2}$. It follows that for the specific case of Gaussian linear observation models, the canonical divergence $D_\phi$ is equivalent to the specific loss identified in \cite{Lucka:2014}. 

\red{We also mention at this point that for Gaussian denoising models; i.e., $\phi(x) = \|y-x\|_2^2 / 2\sigma^2 + \psi(x)$, the estimator $\hat{x}_{D_\phi} = \hat{x}_{\text{MAP}}$ results from the computation of the proximal operator $\textrm{prox}_{\sigma^2 \psi} (y) = \argmin_{x \in \mathbb{R}^n} \|y-x\|_2^2 /2\sigma^2 + \psi(x)$ \cite{PockActaNumerica}. This is equivalent to a gradient step on the Moreau-Yoshida regularisation of $\psi$; i.e., $\hat{x}_{\text{MAP}} = y + \sigma^2 \nabla \tilde{\phi}(y)$, with $\tilde{\psi}(y) = \inf_{x \in \mathbb{R}^n} \|y-x\|_2^2 /2\sigma^2 + \psi(x)$. In like manner, $\hat{x}_{D^*_\phi} = \hat{x}_{\text{MMSE}}$ can be expressed as the gradient step $\hat{x}_{\text{MMSE}} = y + \sigma^2 \nabla \bar{\phi}(y)$, where $\bar{\psi} = \log \int \exp\{-\|y-x\|_2^2 / \sigma^2 - \psi(x) \} \textrm{d}x$ is a different smooth approximation of $\psi$ (please see \cite{Milanfar2018} for details).}

\red{Also note that a different Bregman divergence, namely the KL divergence $KL(u,x) = \int \log\left[\frac{p(y|x)}{p(y|u}\right]p(y|x)\textrm{d}y$, is often used in Bayesian point estimation to define an estimator that is independent of the parametrisation of the likelihood \cite{cprbayes}. This estimator is particularly relevant when the object of interest is $p(y|x)$, as opposed to the value of $x$ itself, for example in prediction problems. This estimator is not often used in imaging sciences.}

{Finally, we notice that because $D_\phi (u,x)$ is derived from $p(x|y)$ may depend on the value of $y$, which is controversial in some lines of Bayesian thinking because it implies that the decision problem underpinning the estimator is defined a-posteriori. \red{This can happen for example in problems involving non-Gaussian observation models.} Our view on this matter is that although decision problems are generally defined a-priori, the case of Bayesian estimators is particular because the decision involved is precisely how to summarise $p(x|y)$, and this decision can be considered a-posteriori if this allows delivering an estimator with favourable accuracy or computational properties. Of course, loss functions that do not depend on the model considered also have advantages, namely the mean squared error loss that also leads to an estimator with good properties (albeit often very expensive to compute). In any case, it is fundamental that one understands how the estimator that one uses summarises $p(x|y)$, and the aim of this work is to improve our understanding of the widely used MAP estimator.}

\subsection{Error bounds for MAP and MMSE estimation}
Theorem \ref{Theo2} establishes that under certain conditions $\hat{x}_{\text{MAP}}$ is a proper Bayesian estimator. Following on from this, it is natural to study the accuracy of $\hat{x}_{\text{MAP}}$ as a Bayesian estimator. The Bayesian approach to this question is to infer the accuracy of $\hat{x}_{\text{MAP}}$ according to the posterior distribution $p(x|y)$. For $\hat{x}_{\text{MMSE}}$ this generally corresponds to computing the expected MSE loss, related to the posterior covariance. This type of analysis can be useful, for example, to identify high dimensional stability conditions (i.e., conditions under which the error grows linearly with $n = \textrm{dim}(x)$). 

Here we perform this type of analysis for $\hat{x}_{\text{MAP}}$ w.r.t. the canonical loss. Precisely, we establish universal estimation error bounds w.r.t. the dual error function $D^*_{\phi}(u,x)$ for MAP and MMSE estimation, for which we have the following result: \vspace{-0.2cm}\\
\begin{proposition}[Expected error bound]\label{Theo3}
Suppose that $\phi(x) = -\log p(x|y)$ is convex on $\mathbb{R}^n$ and $\phi \in \mathcal{C}^1$. Then, 
$$
\textnormal{E}_{x|y}\left[{D^*_\phi(\hat{x}_{\text{MMSE}},x)}\right] \,\leq\, \textnormal{E}_{x|y}\left[{D^*_{\phi}(\hat{x}_{\text{MAP}},x)}\right] \,\leq \, n.
$$
\end{proposition}
\noindent Proof. The proof is reported in the appendix.


\red{We read Proposition \ref{Theo3} as a high dimensional stability result for MAP and MMSE estimation, stating that the expected estimation error, as measured by the dual loss $D^*_{\phi}$, grows at most linearly with the number of image pixels. Therefore, even if the likelihood $p(y|x)$ is poorly identifiable because $\textrm{dim}(y) \ll \textrm{dim}(x)$, or because the linear operator $A$ is very rank deficient, or because $y$ is corrupted by Poison noise, in smooth log-concave settings the expected error cannot grow polynomially or with a linear constant greater than $1$.}

\red{To formally study the expected error as $n$ increases we consider a generic log-concave stochastic process $\mathbb{X} = \{x^{(n)}, n \in \mathbb{N}\}$, where for each $n \in \mathbb{N}$, the random vector $x^{(n)} = (x_1,\ldots,x_n) \in \mathbb{R}^n$ has marginal distribution $p_n(x^{(n)}|y) \propto \exp\{-\phi_n (x^{(n)})\}$ for some convex function $\phi_n : \mathbb{R}^n \rightarrow (-\infty, \infty]$. We also assume that the entropy rate of $\mathbb{X}$ is finite; i.e., $\lim_{n\rightarrow\infty} \textnormal{E}_{x^{(n+1)}|y}[\phi_{n+1}({x}^{(n+1)})] - \textnormal{E}_{x^{(n)}|y}[\phi_n({x}^{(n)})] < \infty$ \cite{CoverBook}. This limit captures the asymptotic information gain per pixel and characterises global statistical features of the image, particularly correlations at any rage. In log-concave settings, this condition holds for example when $\lim_{n\rightarrow\infty} \phi_n(\hat{x}^{(n)}_{\text{MAP}})/n < \infty$ ; it also holds when $\mathbb{X}$ is strongly stationary \cite{Bobkov2011b}. By analysing Proposition \ref{Theo3} in this setting we see that
$$
\textnormal{E}_{x^{(n)}|y}\left[{D^*_{\phi_n}(\hat{x}^{(n)}_{\text{MMSE}},x^{(n)})}\right] \,\leq\, \textnormal{E}_{x^{(n)}|y}\left[{D^*_{\phi_n}(\hat{x}^{(n)}_{\text{MAP}},x^{(n)})}\right] \,\leq \, n \,.
$$ 
Then, because the entropic rate of $\mathbb{X}$ is finite $\lim_{n\rightarrow\infty} \textnormal{E}_{x^{(n)}|y}[\phi_n({x}^{(n)})]/n < \infty$  \cite{CoverBook}, and hence the dimension-normalised expected errors verify
$$
\lim_{n\rightarrow\infty} \textnormal{E}_{x^{(n)}|y}\left[{D^*_{\phi_n}(\hat{x}^{(n)}_{\text{MMSE}},x^{(n)})}/n\right] \,\leq\, \lim_{n\rightarrow\infty} \textnormal{E}_{x^{(n)}|y}\left[{D^*_{\phi_n}(\hat{x}^{(n)}_{\text{MAP}},x^{(n)})}/n\right] \,\leq 1\, .
$$ 
We emphasise again this form of dimension stability is not generally available in estimation problems, and is a direct consequence of the log-concavity of the model and its relationship with the MAP and MMSE estimators. Finally, observe that the above error bounds are tight; e.g., the trivial i.i.d. process $x_i | y \sim \textrm{Exp}(\lambda_y)$, for $i \geq 1$, $\lambda_y \in \mathbb{R}^+$, attains the bound. Lastly, we conjecture that this bound can be improved for specific subclasses of log-concave models by using entropy rate results from the probability literature; future work will investigate this.}

\subsection{Connections to other works}\label{connections}
{We conclude this section by discussing some connections between this paper and other theoretical works related to MAP estimation. As previously explained, Theorem \ref{Theo2} directly builds on \cite{Lucka:2014}, which considered the class of log-concave models $p(x|y) \propto \exp\{-\|y-Ax\|^2_{\Sigma^{-1}}/2 - \lambda h(x)\}$ with Gaussian likelihood $y \sim \mathcal{N}(Ax,\Sigma)$ and prior $p(x) \propto \exp\{\lambda h(x)\}$, and establishes that in this case $\hat{x}_{\text{MAP}}$ is the Bayesian estimator for the Bregman loss $L(u,x) = \|A(u-x)\|^2_{\Sigma^{-1}} + 2\lambda D_h(u,x)\}$. Theorem \ref{Theo2} generalises this result to a larger class of posterior distributions and provides motivation for this unconventional loss function by showing that it stems directly from the consideration of the geometry of the parameter space. Proposition \ref{Theo3} provides further motivation for this loss by establishing explicit bounds on the expected estimation error.}

\red{It is worth mentioning that the generalisation of \cite{Lucka:2014} to other log-concave models was developed simultaneously and independently in Burger et al. \cite{Burger2016} (see \cite[Theorem 3.2]{Burger2016} for MAP estimation, and \cite[Theorem 4.3] {Burger2016} for MSSE estimation, which is also closely related to \cite[Proposition 1]{Banerjee:2005}). Moreover, that work also analyses the expected estimation error involved in MAP and MMSE estimation but w.r.t. other divergence functions. More precisely, Burger et al. \cite{Burger2016} uses the Bregman divergence $D_h$ related to the regulariser or negative log-prior, whereas we use the canonical Bregman divergence $D_\phi$ related to the negative log-posterior. As mentioned previously, $D_{h}$ grasps important features of $x$ (e.g;, sparsity, regularity, smoothness), and is always independent of the observed data $y$, whereas $D_\phi$ is independent of $y$ only in special cases (e.g., Gaussian linear observation models).}

Moreover, Burger et al. also show that $\textnormal{E}_{x|y}\left[{D^*_h(\hat{x}_{\text{MMSE}},x)}\right] \,\leq\, \textnormal{E}_{x|y}\left[{D^*_{h}(\hat{x}_{\text{MAP}},x)}\right]$ and conclude that $\hat{x}_{\text{MMSE}}$ outperforms $\hat{x}_{\text{MAP}}$ when the estimation error is measured in this way, which is independent of $y$. \red{To analyse how these expected errors behave as dimensionality increases we combine this result with Proposition \ref{Theo3} and obtain the following bound:}

\begin{corollary}\label{Theo4}
Suppose that $h(x) = -\log p(x)$ is convex on $\mathbb{R}^n$ and $\phi \in \mathcal{C}^1$. Then, 
$$
\textnormal{E}_{x|y}\left[{D^*_h(\hat{x}_{\text{MMSE}},x)}\right] \,\leq\, \textnormal{E}_{x|y}\left[{D^*_{h}(\hat{x}_{\text{MAP}},x)}\right] \leq n\,.
$$
\end{corollary}
\noindent Proof. The proof follows directly from combining \cite[Theorem 5.1]{Burger2016} with Proposition \ref{Theo3} and the fact that $D^*_{h} \leq D^*_{\phi}$, for any splitting $\phi = h + f$ where $h$ and $f$ are convex functions. \red{This result can also be derived from the integration by parts argument in \cite{Burger2016}.}

\red{Again, we read Corollary \ref{Theo4} as a high dimensional stability result for MAP and MMSE estimation, stating that the expected estimation error, measured in this case by the dual loss $D^*_{h}$, grows at most linearly with the number of image pixels. Polynomial growth or faster linear growth is not possible within the class of smooth log-concave models, even if the problem is very ill-conditioned. At the same time, this linear rate cannot be further improved, as any i.i.d. process $p(x|y) \prod_{i=1}^n p(x_i | y)$ will have an error that grows linearly with $n$.}

\section{Illustrative examples}
As a way of illustrating our theory, we now analyse the geometry of a simple image denoising model in the wavelet domain, and of the regularisation-by-denoising (RED) Bayesian models recently proposed in \cite{RED:2017}.
\subsection{Wavelet image denoising model}
{In this example we analyse the behaviour of MAP estimation in linear inverse problems with sparsity-promoting or shrinkage priors. Without loss of generality, we first consider a simple additive noise observation model $y = x + w$, with noise $w \sim \mathcal{N}(0,\sigma^2 I_n)$ with variance $\sigma^2 \in \mathbb{R}^+$, that allows a detailed analysis. To recover $x$ we put a shrinkage prior on a wavelet representation $z = W x$ of $x$, where $W$ is some orthogonal wavelet transform. More precisely, we use the smoothed Laplace prior 
\begin{equation}\label{smoothedL1}
p(z) \propto \exp\{- \lambda\sum_{i=1}^n \sqrt{z^2_i + \alpha^2}\}\, ,
\end{equation}
where $\lambda \in \mathbb{R}^+$ and $\alpha \in \mathbb{R}^+$ are respectively scale and shape regularisation parameters; this prior is also known as the pseudo-Huber, Hardy, or Charbonnier prior \cite{GEM2006,Charbonnier}. The likelihood is $p(y|z) \propto \exp\{-\frac{1}{2\sigma^2}\|y - W^\top z\|_2^2\}$, and hence the posterior for the wavelet coefficients is
\begin{equation*}
\begin{split}
p(z|y) &\propto \exp\{-\frac{1}{2\sigma^2}\|y - W^\top z\|_2^2 -\lambda\sum_{i=1}^n \sqrt{z^2_i + \alpha^2}\}\, . \\
\end{split}
\end{equation*}
To check that Theorem \ref{Theo2} and Proposition \ref{Theo3} apply, we note that $\phi(z) = -\frac{1}{2\sigma^2}\|y - W^\top z\|_2^2 -\lambda\sum_{i=1}^n \sqrt{z^2_i + \alpha^2}$ belongs to $\mathcal{C}^\infty$ and has a diagonal Hessian matrix with elements given by
$$
\frac{\partial^2}{\partial z^2_i} \phi(z) = \frac{1}{\sigma2} + \frac{\lambda  \alpha^2}{(\alpha^2 + z_i^2)^{3/2}} \, .
$$
Noticing that the elements $\tfrac{\partial^2}{\partial z^2_i} \phi(z)$ take values in $[\frac{1}{\sigma2},\frac{1}{\sigma2}+\lambda]$ for all $z \in \mathbb{R}^n$, we conclude that $\phi(z)$ is strongly convex. Notice that, similarly to the previous example, the geometry of the manifold $\{\mathbb{R}^n,g\}$ does not depend on the value of the observation $y$, and hence the canonical divergences are independent of $y$ too.}

{Because of the action of the shrinkage prior \eqref{smoothedL1}, the Bayesian model $p(z|y)$ will promote solutions that have some large wavelet coefficients and some coefficients close to zero. This behaviour is controlled by the regularisation parameter $\lambda$, and also by the choice of the Bayesian estimator used to summarise $z|y$. In particular, MAP estimation may significantly accentuate shrinkage. This can be theoretically analysed in different ways, and in particular by using Theorem \ref{Theo2}. Accordingly, $\hat{z}_{MAP}$ minimises the expected canonical divergence on $\{\mathbb{R}^n, g\}$, given by the $\phi$-Bregman divergence
\begin{equation*}
\begin{split}
D_\phi (u,z) =& \,\phi(u) - \phi(z) - \nabla \phi(z)^\top (u-z)\, ,\\
= & \,\tfrac{1}{2\sigma^2}\|W^\top u - W^\top z\|^2_2 + \lambda \sum_{i=1}^n \left[\sqrt{u^2_i + \alpha^2} - \sqrt{z^2_i + \alpha^2} + \frac{z^2_i - z_i u_i}{\sqrt{z^2_i + \alpha^2}}\right]\, .\\
\end{split}
\end{equation*}
}
{Because $WW^\top = {I}_n$ we have that $D_\phi$ is fully separable, i.e., $D_\phi (u,z) = \sum_{i=1}^n D_\psi (u_i,z_i)$ with
$$
D_\psi (u_i,z_i) = \tfrac{1}{2\sigma^2} (u_i - z_i)^2 + \lambda \frac{\sqrt{z^2_i + \alpha^2}\sqrt{u^2_i + \alpha^2} - z_i u_i - \alpha^2}{\sqrt{z^2_i + \alpha^2}}\, .
$$
Because $D_\psi$ is a divergence it promotes values of $u_i$ that are close to $z_i$. To develop an intuition for $D_\psi$ it is useful to analyse its behaviour when $z_i$ is small and when it is large relative to $\alpha$.  Observe that $D_\phi$ has a quadratic term related to the likelihood, and a non-quadratic term related to the shrinkage prior.  When $z_i \gg \alpha$ the non-quadratic term vanishes and hence
$$
D_\psi (u_i,z_i) \approx \tfrac{1}{2\sigma^2} (u_i - z_i)^2\, .
$$
As a result, if the observed data is such that the posterior distribution for $z_i | y$ has most of its mass in large values of $z_i$, the MAP estimate for $z_i$ will essentially coincide with the MMSE estimate given by the posterior mean of $z_i | y$. In this case there is no additional shrinkage from the estimator. Conversely, when $z_i \ll \alpha$ the estimator will significantly boost shrinkage. More precisely, when $z_i \ll \alpha$, the non-quadratic term behaves as\begin{equation*}
\begin{split}
D_\psi (u_i,z_i) &\approx \tfrac{1}{2\sigma^2} (u_i - z_i)^2 + \lambda |u_i| \, ,\\
&\approx \tfrac{1}{2\sigma^2} u_i^2 + \lambda |u_i|\, ,
\end{split}
\end{equation*}
for $u_i \gg \alpha$, and for $u_i \ll \alpha$ as 
\begin{equation*}
\begin{split}
D_\psi (u_i,z_i) &\approx \tfrac{1}{2\sigma^2} (u_i - z_i)^2 + \lambda \left[\frac{u_i^2}{2\alpha} + \frac{z_i^2}{2\alpha} - \frac{- z_i u_i}{\alpha}\right] = \left(\frac{1}{2\sigma^2} + \frac{\lambda}{2\alpha}\right)(u_i - z_i)^2_2 \, , \\
& \approx \left(\frac{1}{2\sigma^2} + \frac{\lambda}{2\alpha}\right) u_i^2\, .
\end{split}
\end{equation*}
In these two cases $D_\psi$ strongly promotes $u_i$ values that are close to zero,  either explicitly via the shrinkage term $\lambda |u_i|$, or by amplifying the convexity constant of the quadratic loss from $1/{2\sigma^2}$ to $1/{2\sigma^2} + {\lambda}/{2\alpha}$. As a result, if the posterior distribution for $z_i | y$ has mass in small values of $z_i$, then the MAP estimate will intensify the shrinkage effect of the prior. This additional shrinkage is not observed with other loss functions, e.g., MMSE, and is consistent with the empirical observation that MAP estimation performs well with shrinkage priors.
}

{For illustration, Figure \ref{wavelet} shows an experiment with the \texttt{Flinstones} image of size $256 \times 256$ pixels. Figure \ref{wavelet}(a) shows a corrupted observation $y = x + w$ with noise $w \sim \mathcal{N}(0,\sigma)$ with $\sigma = 0.08$, which has a signal-to-noise ratio of $17.6$dB. The restored imaged obtained by MAP estimation is displayed in Figure \ref{wavelet}(b), this estimate has a signal-to-noise ratio of $19.8$dB (we used a Haar wavelet decomposition with four levels and $\lambda = 12$ and $\alpha = 0.01$ for all scales except the coarse scale, for which we used a Jeffreys' prior $p(z_i) \propto 1$ to avoid excessively biasing the estimates). For comparison, we also report $\hat{x}_{MMSE}$, which in this experiment performs poorly (signal-to-noise ratio of $17.7$dB). The MAP estimator obtained with a conventional Laplace or $\ell_1$ prior (i.e., with $\alpha \rightarrow 0$) has a worse signal-to-noise ratio ($18.8$dB, not displayed)}

{Because $p(z|y)$ is fully separable, i.e., $p(z|y) \prod_{i=1}^n p(z_i|y)$ - and thus $D_\psi = \sum_{i=1}^n D_\psi (u_i, z_i)$ - the action of these estimator can be clearly visualised by plotting the estimation function that performs the denoising of the wavelet coefficients: for MAP estimation this is given by $\hat{z}_{i,MAP} (y): y \rightarrow  \argmin_{u_i} \textrm{E}_{z_i} \left[D_\psi (u_i, z_i) | w_i^\top y \right]$, where $w_i^\top y$ is the $i$th wavelet coefficient of $y$; and for MMSE estimation it is the marginal posterior mean $\hat{z}_{i,MMSE} (y): y \rightarrow  \int z_i p(z_i|y) \textrm{d}z_i$. These functions, displayed in Figure \ref{wavelet}(d), clearly show the importance of the choice of the loss used to summarise $p(z|y)$.}
{
\begin{figure}\label{wavelet}
	\begin{tabular}{cc}
\includegraphics[width=7.3cm]{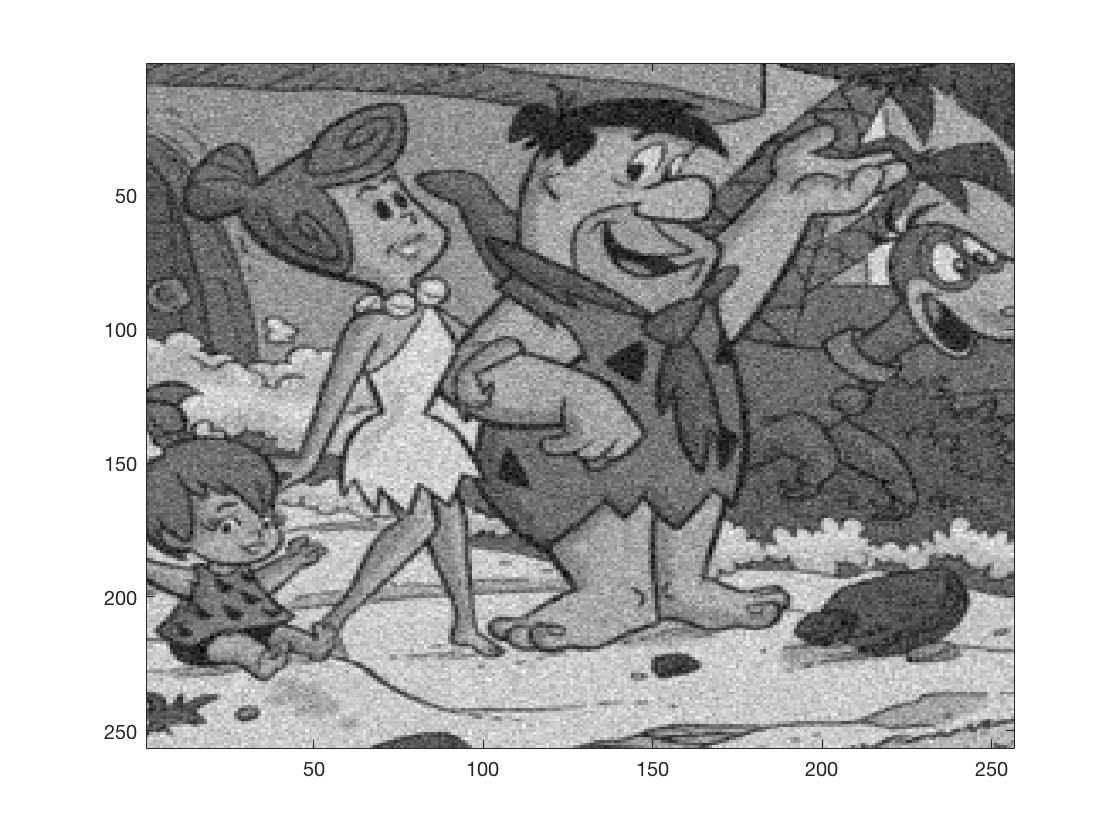}&\includegraphics[width=7.3cm]{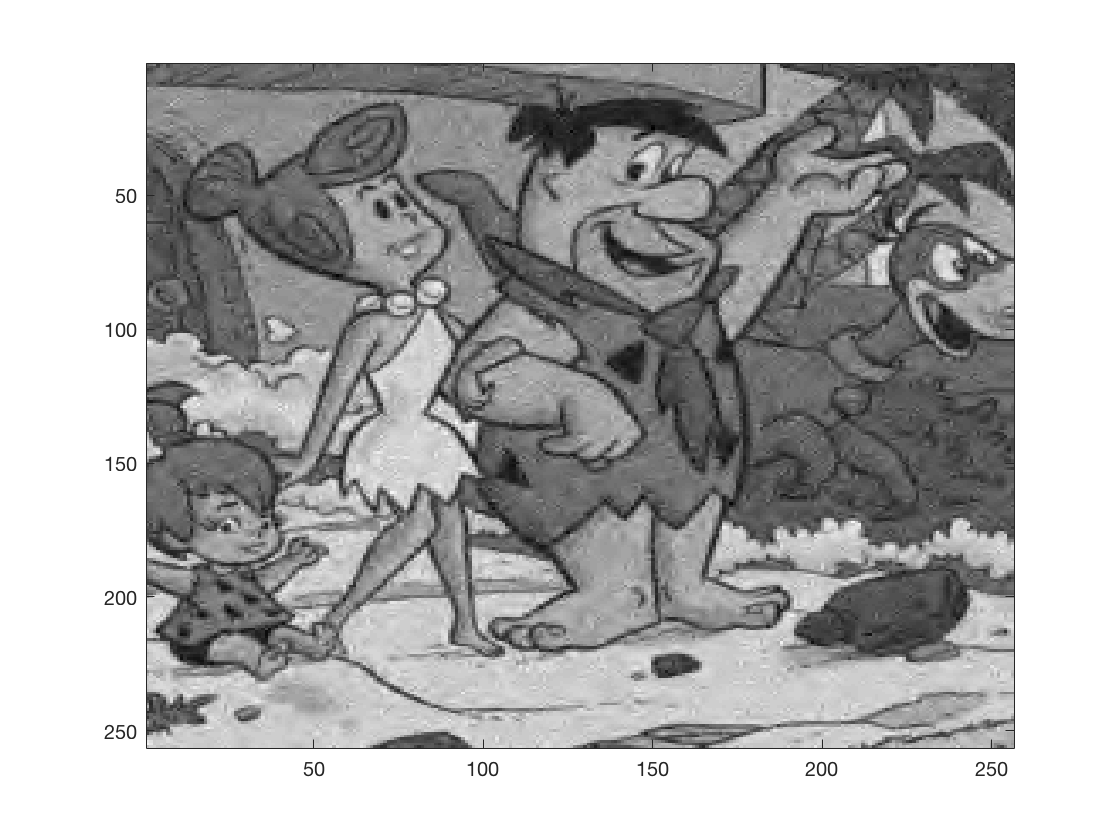}
\\
		{\small (a)  observation $y$} & {\small (b) restored image $\hat{x}_{MAP}$} 
		\\
\includegraphics[width=7.3cm]{wavelet_observation.png}&\includegraphics[width=7.3cm]{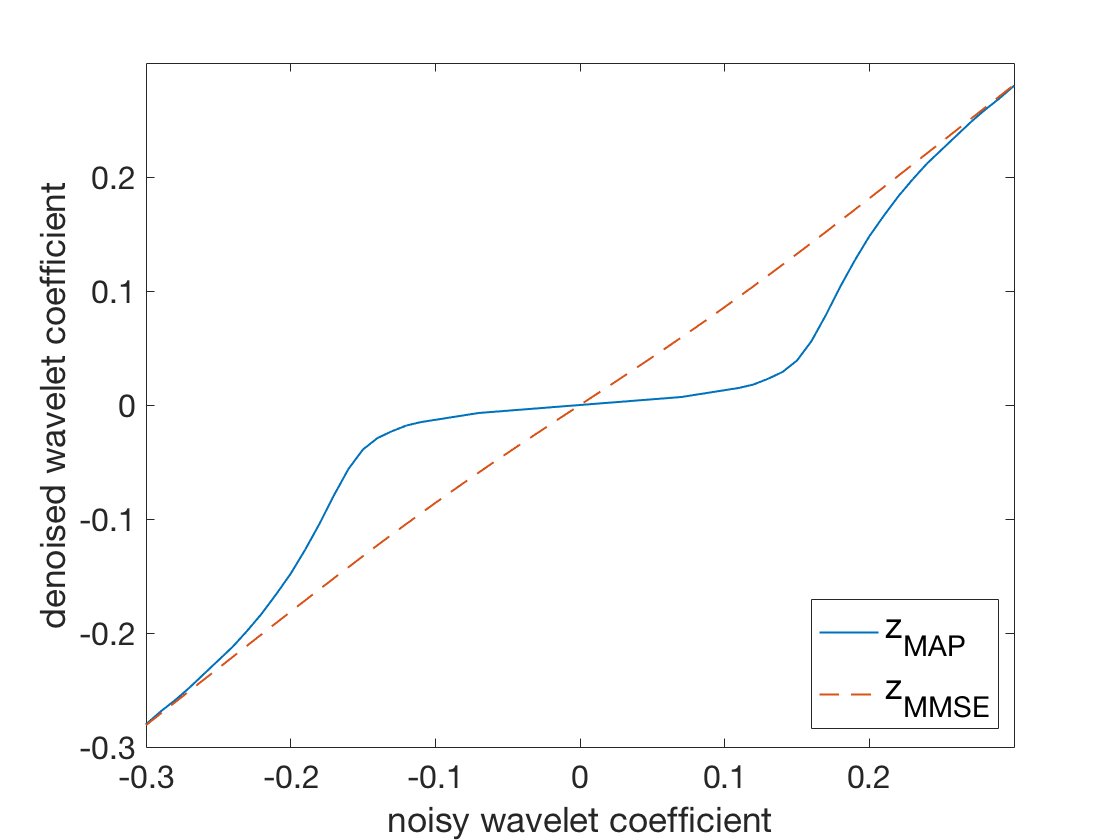}
\\
		{\small (c) restored image $\hat{x}_{MMSE}$} & {\small (d) denoising functions for $\hat{z}_{MAP}$ and $\hat{z}_{MMSE}$}		
		\\
        \end{tabular}
\caption{Wavelet denoising experiment with the \texttt{Flinstones} image with the smoothed Laplace prior \eqref{smoothedL1}.}
\end{figure}
}

{We emphasise at this point that this experiment has been selected to highlight the additional shrinkage obtained by using MAP estimation instead of MMSE estimation. However, there are other models where, because of the likelihood and the choice of the wavelet representation and the parameters used, shrinkage arises mainly from the prior. In that case MAP and MMSE estimation perform equally well. To illustrate this point, Figure \ref{wavelet2} below shows the reconstruction results obtained in \cite{Cai2018} with MAP and MMSE estimation for a radio-interferometric imaging problem with a very similar model of the form $p(z|y) \propto \exp\{-\frac{1}{2\sigma^2}\|y -  A W^\top z\|_2^2 -\lambda\sum_{i=1}^n |z_i |\}$, where the difference is that the likelihood term involves a linear operator $A$ modelling the radio-telescope system (see \cite{Cai2018} for more details about the model and the algorithms used to compute the estimates). Observe that in this case both MAP and MMSE estimation deliver excellent and remarkably similar results. A similar empirical observation is reported in \cite{Lucka:2014} for a sparse tomography experiment using the Besov wavelet model of \cite{Siltanen}, which is closely related to the model considered here.
}

\begin{figure}\label{wavelet2}
	\begin{tabular}{cc}
\includegraphics[width=7.3cm]{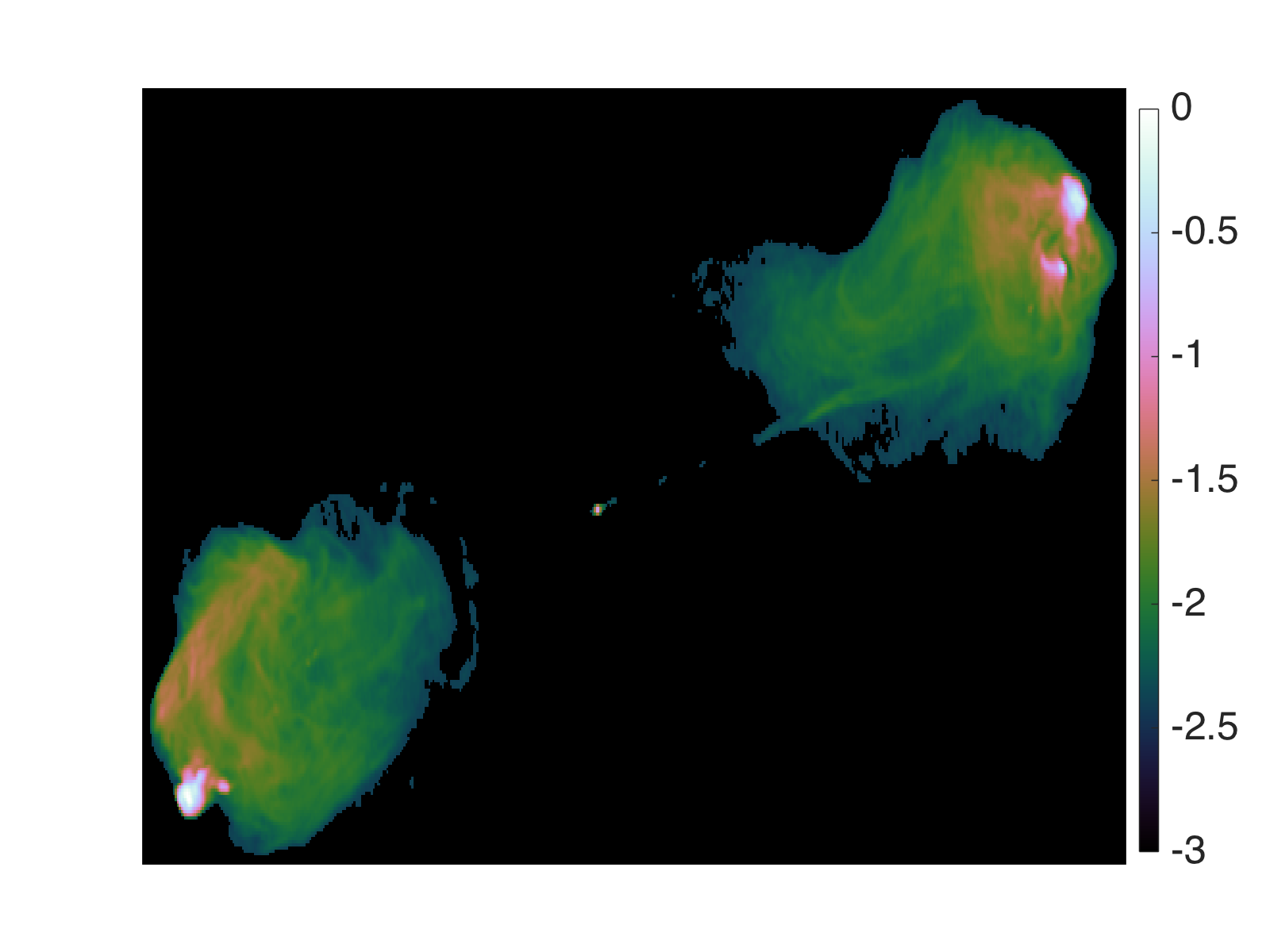}&\includegraphics[width=7.3cm]{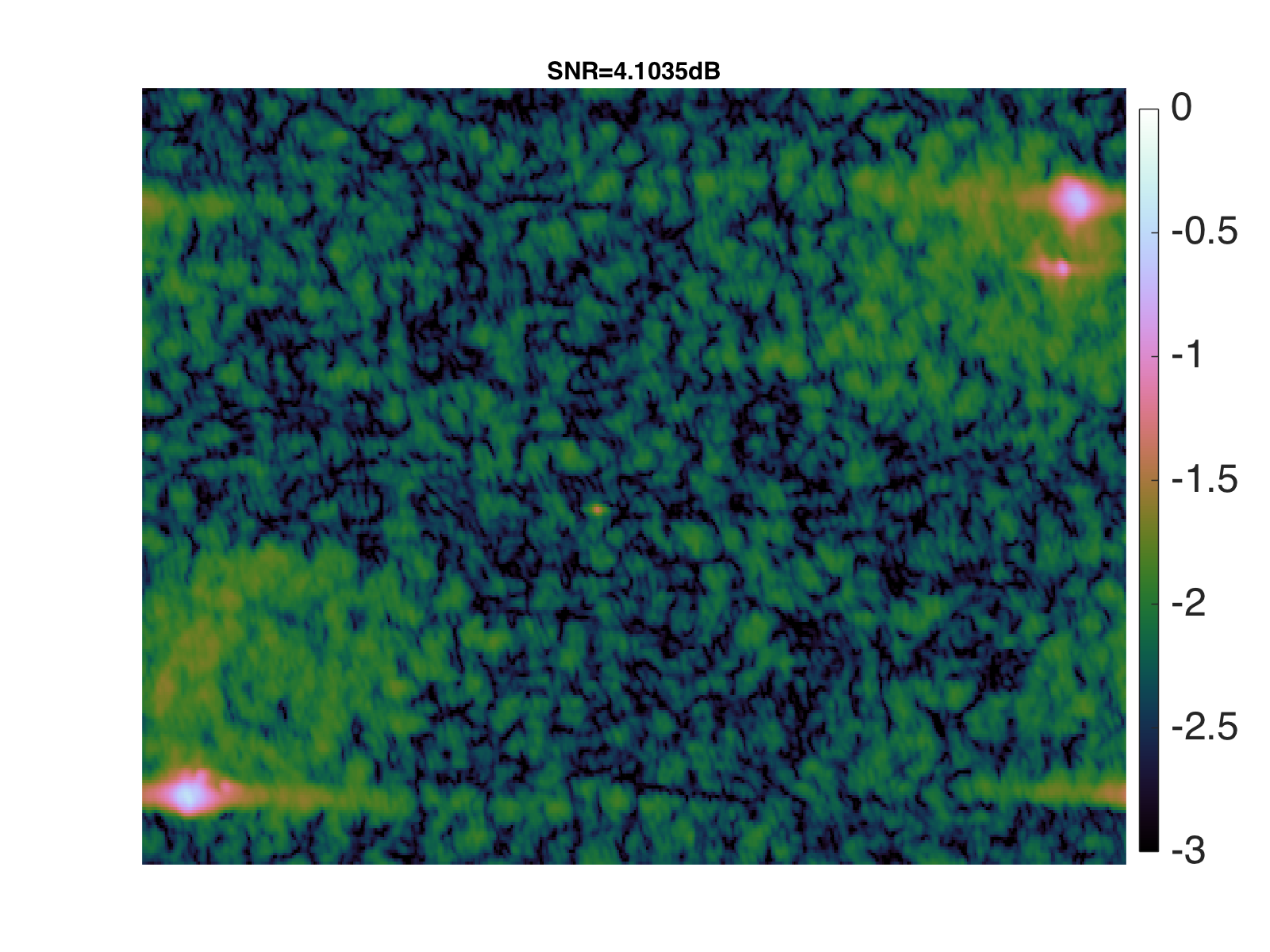}
\\
		{\small (a)  true image $x$} & {\small (b) observation $y$} 
		\\
\includegraphics[width=7.3cm]{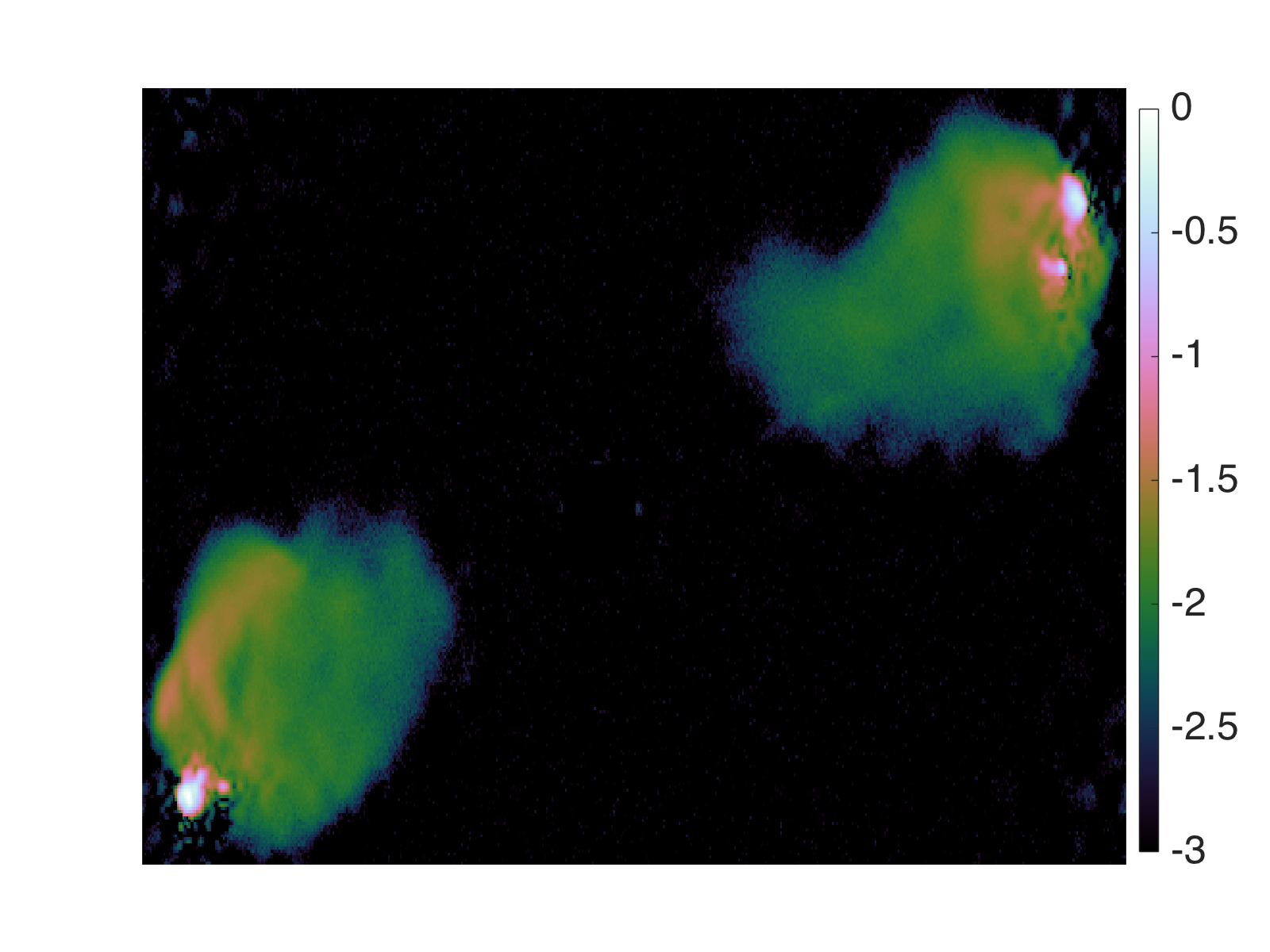}&\includegraphics[width=7.3cm]{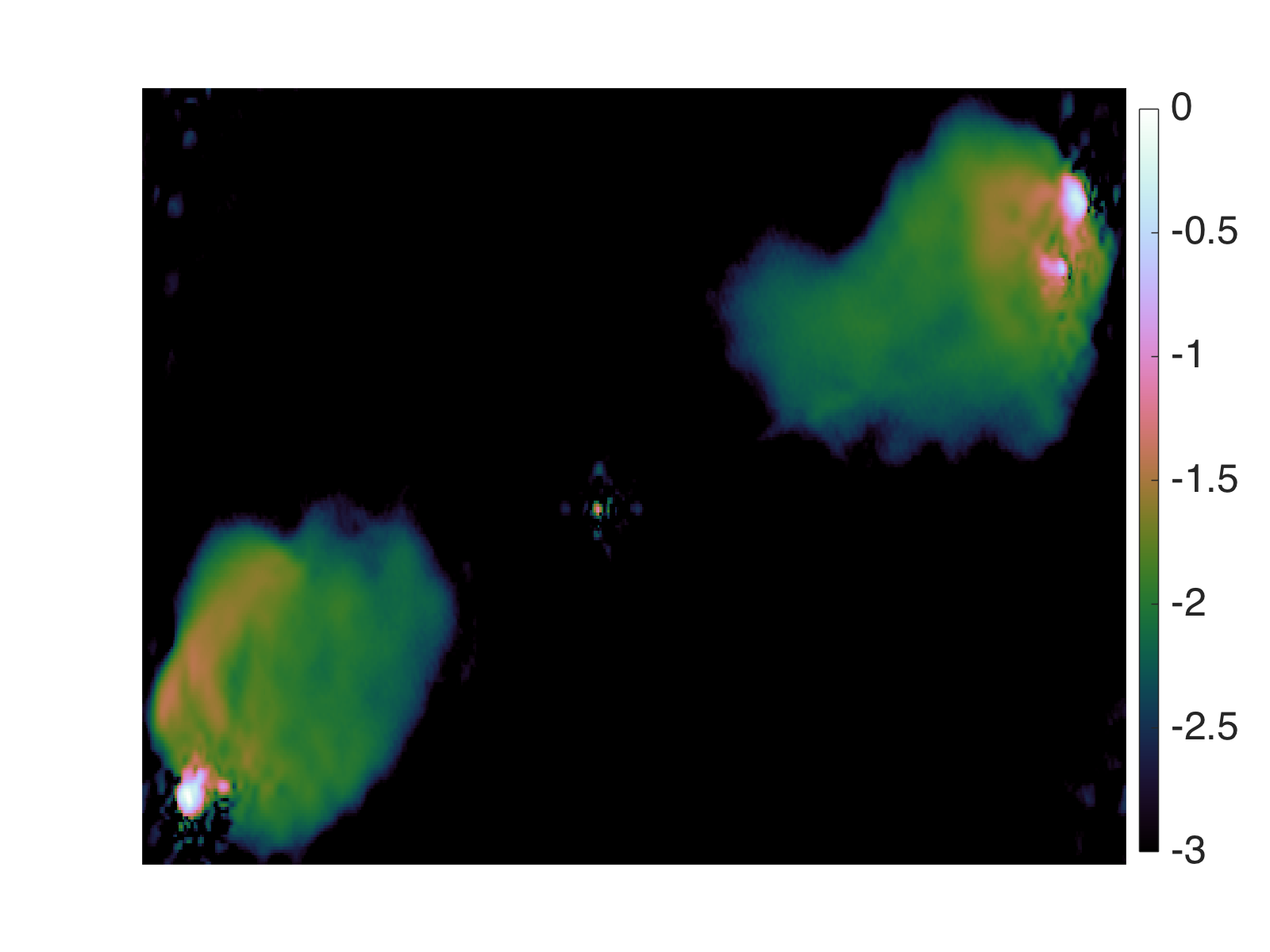}
\\
		{\small (c) restored image $\hat{x}_{MAP}$} & {\small (d) restored image $\hat{x}_{MMSE}$}		
		\\
        \end{tabular}
\caption{Bayesian radio-interferometric imaging experiment with the \texttt{Cygnus A} radio galaxy (size $256 \times 512$ pixels). See \cite{Cai2018} for more details.}
\end{figure}

{To conclude, shrinkage priors promote solutions that are sparse or approximately sparse through two mechanisms: directly through the definition of the bayesian model $p(z|y)$, and indirectly through the loss function used to summarise $z|y$. In the case of MAP estimation, this loss function is a Bregman divergence that can significantly amplify shrinkage. In some cases this may lead to better estimation performance. Therefore, when designing Bayesian procedures for imaging problems, it is important to carefully consider both the model and the Bayesian estimator used to summarise the information provided by the model.}

\subsection{Regularisation-by-denoising (RED) Bayesian models}
As a second illustrative example we analyse the geometry of the RED Bayesian models \cite{RED:2017}. In this class of models the prior $p(x)$ is defined implicitly through an image denoising algorithm. {Precisely, starting from some image denoising filter $f : \mathbb{R}^n \rightarrow \mathbb{R}^n$, we posit the prior}
\begin{equation}\label{REDprior0}
p(x) \propto \exp\{-\tfrac{\lambda}{2} x^\top[x-f(x)]\}\, ,
\end{equation}
{which promotes values of $x$ that are approximately invariant to filtering by $f$ (i.e., for which $f(x) \approx x$), the rationale being that these are values that $f$ considers to be realistic images.} Note that this approach has close connections to {plug-and-play} priors, that are also defined through denoising algorithms \cite{PPP:2015}.  

{The RED framework \cite{RED:2017} assumes that $f$ verifies three regularity conditions that are necessary to make inference with \eqref{REDprior0} analytically and computationally tractable. First, $f$ is smooth, at least $\mathcal{C}^2$. Second,  the Jacobian matrix $J_f(x)$ is symmetric and has all its eigenvalues in $[-1,1]$ for all $x \in \mathbb{R}^n$. Third, $f$ is locally homogenous, i.e., $\lim_{\epsilon \rightarrow 0} \epsilon^{-1} [f(x + \varepsilon x) - f(x)] = f(x)$ $\forall x \in \mathbb{R}^n$; this property implies that $f(x) = J_f(x) x$. Under these assumptions, it is possible to express \eqref{REDprior0} in the following pseudo-quadratic form
\begin{equation}\label{REDprior}
p(x) \propto \exp\{-\tfrac{\lambda}{2} x^\top \Lambda_f(x) x\}\,.
\end{equation}
where $\Lambda_f(x) = {I}_n - J_f(x)$ plays the role of an image-adapted graph Laplacian operator, highlighting the connection between $f$ and the model geometry \cite{RED:2017}.}

{Notice that these regularity assumptions imply the regulariser or negative log-prior $h(x) = -\log p(x)$ is $\mathcal{C}^3$ and convex, which is important for example for the efficient computation of $x_{MAP}$ by optimisation. To show that $h(x) \in \mathcal{C}^3$ we first use the fact that $\nabla h(x) = \nabla \log p(x) = -\lambda[x - f(x)/2 - J_f(x)^\top x/2] = -\lambda[x-f(x)]$, where we have used the symmetry $J_f(x)^\top = J_f(x)$ and the identity $f(x) = J_f(x) x$ related to the local homogeneity of $f$. Therefore, the Hessian matrix of $h(x)$ has elements given by $\partial_i \partial_j h(x) = - \Lambda_f(x)_{i,j} = J_f(x)-I_n$, which are also continuously differentiable because $f \in \mathcal{C}^2$, and hence $h(x) \in \mathcal{C}^3$. Moreover, the Hessian matrix of $h(x)$, given by $-\Lambda_f(x) = J_f(x)-I_n$, is negative-semidefinite because $J_f(x)$ has all its eigenvalues in $[-1,1]$ $\forall x \in \mathbb{R}^n$, and consequently $h(x)$ is convex. As a result, if the negative log-likelihood $-\log p(y|x)$ is $\mathcal{C}^3$ and convex w.r.t. $x$, then $\phi(x) = -\log p(x|y)$ is also $\mathcal{C}^3$ and log-concave, and Theorem \ref{Theo2} and Proposition \ref{Theo3} apply.} 

As illustrative example, consider linear inverse problems of the form $y = Ax + w$, where $A$ is a known linear operator, and $w \sim \mathcal{N}(0,\sigma^2 I_n)$ with noise variance $\sigma^2 \in \mathbb{R}^+$. The resulting RED Bayesian model has posterior density\footnote{Because $\Lambda(x)$ is potentially rank deficient, to guarantee that $p(x | y)$ is a proper probability density function we further assume that $\textrm{ker}(A^\top A) \cap \textrm{ker}(\Lambda(x)) = \{0\}$ for all $x \in \mathbb{R}^n$.}
\begin{equation}\label{REDmodel}
p(x|y) \propto \exp\{- \tfrac{1}{2\sigma^2}\|y-Ax\|^2_2 - \tfrac{\lambda}{2} x^\top \Lambda(x) x\}\, .
\end{equation}
This distribution is strongly log-concave and $\mathcal{C}^3$, and hence Theorem \ref{Theo2} and Proposition \ref{Theo3} apply. More precisely, we have $\phi(x) =
\tfrac{1}{2\sigma^2}\|y-Ax\|^2_2 + \tfrac{\lambda}{2} x^\top \Lambda(x) x$, which induces the metric
\begin{equation}\label{hessian_REDmodel}
g_{i,j}(x) = \partial_i \partial_j \phi (x) = \{\sigma^{-2} A^\top A  + \lambda \Lambda(x)\}_{i,j}\, .
\end{equation}
Observe that \eqref{hessian_REDmodel} combines an Euclidean geometry term $A^\top A$ from the Gaussian likelihood, and a non-Euclidean term from the Laplacian $\Lambda(x)$. Again, note that for this class of models the geometry of the manifold $\{\mathbb{R}^n,g\}$ does not depend on the value of the observation $y$. 

{Moreover, from Theorem \ref{Theo2}, the estimator $\hat{x}_{\text{MAP}}$ is the Bayesian estimator associated with the canonical divergence on $\{\mathbb{R}^n, g\}$, given by the $\phi$-Bregman divergence
\begin{equation*}
\begin{split}
D_\phi (u,x) =& \,\phi(u) - \phi(x) - \nabla \phi(x)^\top (u-x)\, ,\\
= & \,\tfrac{1}{2\sigma^2}\|y-Au\|^2_2 + \tfrac{\lambda}{2} u^\top \Lambda(u) u - \tfrac{1}{2\sigma^2}\|y-Ax\|^2_2 + \tfrac{\lambda}{2} x^\top \Lambda(x) x\\
&\,- [\sigma^{-2}A^\top(Ax-y) + \lambda x - \lambda f(x)]^\top (u-x)\, ,\\
=&\, \sigma^{-2} [u^\top A^\top A u + x^\top A^\top A x - 2 u^\top A^\top A x] + \lambda [u^\top \Lambda(u) u + x^\top \Lambda(x) x - 2 u^\top \Lambda(x) x]\, ,\\
=&\, \sigma^{-2} D_{A^\top A}(u,x) + \lambda D_{\Lambda} (u,x)\, ,
\end{split}
\end{equation*}
where the Mahalanobis (Euclidean) distance
\begin{equation*}
\begin{split}
D_{A^\top A} (u,x) &= \|u-x\|_{A^\top A}\, , \\
&= u^\top A^\top A u + x^\top A^\top A x - 2 u^\top A^\top A x\, ,
\end{split}
\end{equation*}
is a measure of prediction MSE related to the Gaussian likelihood, and
\begin{equation*}
\begin{split}
D_{\Lambda} (u,x) = u^\top \Lambda(u) u + x^\top \Lambda(x) x - 2 u^\top \Lambda(x) x\, ,
\end{split}
\end{equation*}
is related to the Laplacian $\Lambda(x)$, which encodes the geometry of the prior (to compute $D_{\Lambda}$ we used the local homogeneity property $f(x) = J_f(x) x$ of the denoiser, see \cite{RED:2017} for details)}.

Finally, observe that $D_{\Lambda}$ is very similar to the Euclidean loss $D_{A^\top A}$ in that it measures the difference between the squared norms of $u$ and $x$ and the projection of $u$ on $x$, with the only difference being that for $D_{\Lambda}(u,x)$ these norms and projections are measured on the tangent spaces $\mathcal{T}_u \mathbb{R}^n$ and $\mathcal{T}_x \mathbb{R}^n$ of the manifold $\{\mathbb{R}^n, \Lambda\}$, with inner products specified by $\Lambda$. 

\section{Relaxation of regularity conditions}\label{relax}
We now examine the effect of relaxing the regularity assumptions of Theorem \ref{Theo2}. We consider three main cases: lack of smoothness, lack of strong convexity, and lack of continuity.

\subsection{Non-smooth models}
The results of Theorem \ref{Theo2} hold for non-smooth models with the following modifications. 
{First, assume that $\phi$ is \emph{almost everywhere} $\mathcal{C}^3$ on $\mathbb{R}^n$; i.e., the set of points of $\mathbb{R}^n$ where $\phi$ is not smooth has dimension $n-1$ and hence zero Lebesgue measure. To check that $\phi$ is $\mathcal{C}^3$ almost everywhere it is necessary to analyse the regularity of the second order derivatives $\partial_i \partial_j \phi(x)$ (e.g., if the second derivatives are Lipchitz continuous then $\phi$ is \emph{almost everywhere} $\mathcal{C}^3$ by Rademacher's theorem \cite{Niculescu2018}).} Because in models that are almost everywhere smooth the set of non-differentiable points has no probability mass, this set can be simply omitted in the computation of expectations. Second, because the non-differentiable points do not have Euclidean tangent spaces, instead of a global manifold we need to consider the collection local manifolds associated with the regions of $\mathbb{R}^n$ where $p(x|y)$ is $\mathcal{C}^3$. Each one of these regions has a local canonical divergence given by the Bregman divergence $D(u,x) = D_\phi(u,x) = \phi(u) - \phi(x) - \nabla \phi(x)^\top (u-x)$. Therefore, for these models we need to posit $D_\phi(u,x)$ as the global loss function for any $(u,x) \in \mathbb{R}^n \times  \mathbb{R}^n$ [technically the global loss is the generalised Bregman divergence $D_\phi(u,x) = \phi(u)-\phi(x) - q_x^\top(u-x)$, where $q_x$ belongs to the subdifferential set of $\phi$ at $x$ \cite{CombettesBook}, however the expectation $\textnormal{E}_{x|y}[D_\phi(u,x)]$ is taken over the points where $\phi$ is $\mathcal{C}^3$ and hence $q_x = \nabla \phi(x)$]. By calculating the primal and dual Bayesian estimators related to this global loss we obtain that $\hat{x}_{\text{MAP}} = \argmin_{u \in \mathbb{R}^n} \textnormal{E}_{x|y}[D_\phi(u,x)]$ and $\hat{x}_{\text{MMSE}} = \argmin_{u \in \mathbb{R}^n} \textnormal{E}_{x|y}[D^*_\phi(u,x)]$, similarly to Theorem \ref{Theo2}. Observe that these modifications do not affect the fact that $\hat{x}_{\text{MAP}}$ and $\hat{x}_{\text{MMSE}}$ can correspond to non-differentiable points. Also note that despite not being a global canonical divergence, $D_\phi(u,x)$ is still consistent with the space's Riemannian geometry which is local.

{Many imaging models involve non-smooth norms such as the $\ell_1$ and the nuclear norm, or the total-variation pseudo-norm, that are almost everywhere $\mathcal{C}^1$ but not $\mathcal{C}^3$. More generally, all Lipchitz continuous functions are almost everywhere $\mathcal{C}^1$. In this case, only the second and third parts of Theorem \ref{Theo2} hold. That is, we posit $D_\phi(u,x)$ as the loss function for any $u \in \mathbb{R}^n$ and any $x \in  \mathbb{R}^n$, excluding non-differentiable points, and obtain that $\hat{x}_{\text{MAP}} = \argmin_{u \in \mathbb{R}^n} \textnormal{E}_{x|y}[D_\phi(u,x)]$ and $\hat{x}_{\text{MMSE}} = \argmin_{u \in \mathbb{R}^n} \textnormal{E}_{x|y}[D^*_\phi(u,x)]$ by removing non-differentiable points from the calculation of the expectations. To recover the differential geometric derivation of $D_\phi$ it is necessary to use a smooth approximation of the model, i.e., the smoothed L1 norm $|s| = \sqrt{s^2 + \alpha^2}$ for some arbitrarily small $\alpha > 0$.}

{Finally, also note that the bound $\textnormal{E}_{x|y}\left[{D^*_{\phi}(\hat{x}_{\text{MAP}},x)}\right] \,\leq \, n$ in Proposition \ref{Theo3} is straightforwardly extended to non-smooth models by using the generalised dual Bregman divergence $D^*_{\phi}(u,x) = \phi(x)-\phi(u) - q_u^\top(x-u)$ with subgradient $q_{\hat{x}_{\text{MAP}}} = 0$. Conversely, the other bound $\textnormal{E}_{x|y}\left[{D^*_\phi(\hat{x}_{\text{MMSE}},x)}\right] \,\leq\, n$ is lost (see the appendix for details).}

\subsection{Strictly log-concave models}
For models that are strictly but not strongly log-concave only the second and third results of Theorem \ref{Theo2} remain true. {It is easy to check that the Bayesian estimators w.r.t. $D_\phi = \phi(u) - \phi(x) - \nabla \phi(x)^\top (u-x)$ are still $\hat{x}_{\text{MAP}} = \argmin_{u \in \mathbb{R}^n} \textnormal{E}_{x|y}[D_\phi(u,x)]$ and $\hat{x}_{\text{MMSE}} = \argmin_{u \in \mathbb{R}^n} \textnormal{E}_{x|y}[D^*_\phi(u,x)]$, similarly to strongly log-concave models (see in the appendix that strong log-concavity is not required to prove the second and third parts of Theorem \ref{Theo2})}. Thus, the decision-theoretic derivation of $\hat{x}_{\text{MAP}}$ remains valid, and $\hat{x}_{\text{MAP}}$ and $\hat{x}_{\text{MMSE}}$ remain dual to each other. The high dimensional performance guarantees of Proposition \ref{Theo3} also hold because $\phi$ is convex. However, without strong convexity, $g$ becomes semi-positive definite and $(\mathbb{R}^n,g)$ becomes a singular manifold. Currently, the validity of the interpretation of $D_\phi$ as a canonical divergence in singular manifolds is not clear. The generalisation of canonical divergences and of Theorem \ref{Theo2} to singular manifolds is currently under investigation. {In any case, without strong convexity $D_\phi$ is no longer a divergence in the sense of Definition \ref{divergences} as $D_\phi(x,u)=0$ does not imply $x=u$, which is an important desired property for loss functions.}

\subsection{Models involving constraints}
Finally, in cases where $x|y$ is constrained to a convex set $\mathcal{S} \subset \mathbb{R}^n$ only the first and the third results of Theorem \ref{Theo2} hold. Proceeding similarly to the proof of Theorem \ref{Theo2}, it is easy to show that $D_\phi$ is the canonical divergence of the manifold $(\mathcal{S},g)$, and that the Bayesian estimator related to the dual divergence is $\hat{x}_{\text{MMSE}} = \argmin_{u \in \mathcal{S}} \textnormal{E}_{x|y}[D^*_\phi(u,x)]$. However, the Bayesian estimator that minimises the canonical divergence is now a shifted or tilted MAP estimator
$$
\hat{x}_{D_\phi} = \argmin_{u \in \mathcal{S}} D_\phi(u,\hat{x}_{\text{MAP}}) + u^\top \textnormal{E}_{x|y}[\nabla \phi(x)],
$$
where generally $\textnormal{E}_{x|y}[\nabla \phi(x)] \neq 0$ (see the appendix for details). It is not clear at this point under what conditions $\hat{x}_{\text{MAP}} \approx \hat{x}_{D_\phi}$. Nevertheless, the high dimensional guarantees of Proposition \ref{Theo3} still hold for $\hat{x}_{\text{MAP}}$, providing some theoretical justification for using this estimator.

\subsection{Models with heavy-tails}
{We conclude this section by discussing the difficulties of extending our results to models that are heavy-tailed and hence not log-concave, such as imaging models involving heavy-tailed priors related to compressible distributions \cite{Gribonval2012}. Unfortunately, extending our results to heavy-tailed settings is extremely challenging for several reasons. First, the Hessian matrix of $\phi$ does not define a Riemannian metric because there are regions of the space where it has negative eigenvalues. Also, directly postulating $D_\phi = \phi(u) - \phi(x) - \nabla \phi(x)^\top (u-x)$ as loss function is not appropriate either because $D_\phi$ can take negative values. The analysis is further complicated by the fact that $p(x|y)$ may have an infinite number of maximisers in disconnected areas of the parameter space. As mentioned previously, the derivation of MAP estimation as an approximation arising from the degenerate loss $L_\epsilon(u,x) = \boldsymbol{1}_{\|x-u\|<\epsilon}$ with $\epsilon \rightarrow 0$ also fails in this case \cite{Bassett2016}. Interestingly, MMSE estimation may also struggle here given that models in this class may not have a posterior mean \cite{cprbayes}.}

\section{Conclusion}\label{sec:Conclusion}
MAP estimation is one of the the most successful Bayesian estimation methodologies in imaging science, with a track record of accurate results across a wide range of challenging imaging problems. Our aim here has been to contribute to the theoretical understanding of this widely used methodology, particularly by placing it in the Bayesian decision theory framework that underpins the core Bayesian inference methodologies.

In order to analyse MAP estimators we have adopted an entirely new approach: we allowed the model to specify the loss function, or equivalently the Bayesian estimator, that is used to summarise the information that the model represents. This was achieved by using the connections between model log-concavity, Riemannian geometry, and divergence functions. We first established that if $p(x|y)$ is strongly log-concave, continuous, and $\mathcal{C}^3$ on $\mathbb{R}^n$, then $\phi (x) = -\log p(x|y)$ induces a dually-flat Riemannian structure on the parameter space, where the canonical divergence is the Bregman divergence associated with $\phi$, and where the MAP estimator is the unique Bayesian estimator w.r.t. to this loss function. We also established that the MMSE estimator is the Bayesian estimator w.r.t. the dual canonical loss, and that both estimators enjoy favourable stability properties in high dimensions. We then examined the effect of relaxing these assumptions to models with weaker regularity conditions.

The theoretical results presented in this work provide several valuable new insights into MAP and MMSE estimation. In particular, both estimators stem from Bayesian decision theory and from the consideration of the geometry of the parameter space, and exhibit an interesting form of duality. Also, the expected estimation error - as measured by the canonical loss - is stable in high dimensions; this is in agreement with the remarkable empirical performance observed imaging and other large scale settings. The fact that MAP estimators are available as solutions to convex problems is a fundamental practical advantage. However, our results also show that the predominant view of MAP estimators as hastily approximate inferences, motivated only by computational efficiency, is fundamentally incorrect. We hope that these results will provide some clarity to imaging scientists using MAP estimators, and that they stimulate further research into the theory of this powerful Bayesian methodology.

\section{Acknowledgements}
Part of this work was conducted when the author held a Marie Curie Intra-European Research Fellowship for Career Development at the University of Bristol, and part when he was a visiting professor at the Institut Henri Poincar\'e in France. He is grateful to Yoann Altmann,  Gavin Gibson, Peter Green,  Abderrahim Halimi, Bernd Schroers, Jonty Rougier, and Ben Powell for useful discussion.

\bibliographystyle{siamplain}
\bibliography{refs}

\appendix

\section*{Appendix - Proofs of Theorem \ref{Theo2} and Proposition \ref{Theo3}}\label{sec:Proof}
\subsection*{Proof of Theorem \ref{Theo2}}
The first part of Theorem \ref{Theo2} follows directly from differential geometry and from the regularity properties of $\phi$ (see \cite{amaribook} for an introduction to differential geometry). From differential geometry, under the conditions of Theorem \ref{Theo2}, $\phi$ induces a Riemannian metric on $\mathbb{R}^n$ with coefficients
$$
g_{i,j}(x) = \partial_i \partial_j \phi(x),
$$ 
and where we note that $g(x)$ is positive definite from the strong convexity of $\phi$. Similarly, we have the affine connection coefficients $\Gamma_{ij,\,k} = \partial_i \partial_j \partial_k \phi(x)$.

Moreover, because $\phi$ is convex it endows $(\mathbb{R}^n, g)$ with a dual affine coordinate system $\eta$, related to the primal coordinate system by the duality $\eta_x = \nabla \phi(x)$ and $x_\eta = \nabla \phi^\star(\eta)$, where $\phi^\star(\eta) = \max_{x\in\mathbb{R}^n} x^\top\eta - \phi(x)$ is the convex conjugate of $\phi$  \cite[Ch. 3]{amaribook}. As a result we have a dual Riemannian metric $g^\star$ w.r.t. $\eta$, with coefficients given by
$$
g^\star_{i,j}(\eta) = \partial_i \partial_j \phi^\star(\eta),
$$
and a dual affine connection $\Gamma^\star$ with coefficients given by
$$
\Gamma^\star_{ij,\,k} (\eta)= \partial_i \partial_j \partial_k \phi^\star(\eta).
$$

Finally, it is easy the check that $x$ and $\eta$ are mutually dual w.r.t. $g$. That is, for all $x\in \mathbb{R}^n$
$$
g^\star(\eta_x) = g(x)^{-1}
$$
which implies that $(\mathbb{R}^n, g, \Gamma, \Gamma^\star)$ is a dually-flat Riemannian manifold  \cite[Ch. 3]{amaribook}. {Please see \cite[Section 2]{Nielsen2018} for an excellent introduction to dually-flat structures and their main properties.}

From  \cite{amari:2015}, in such manifolds the $\Gamma$-geodesic connecting $u\rightarrow x$ in \eqref{canonicalDiv} is given by $\gamma_t = u + t(x-u)$, and $\dot{\gamma}_t = x-u$. The proof is then concluded by integration by parts of \eqref{canonicalDiv} to obtain the Bregman divergence $D_\phi (u,x) = \phi(u) - \phi^\star(\eta_x) - \eta^\top_x u$, which also admits the more familiar expression $D_\phi (u,x) = \phi(u) - \phi(x) - \nabla \phi(x)(u - x)$.

To prove the second part of Theorem \ref{Theo2} we use the linearity property of the expectation operator to express the definition $\hat{x}_{D_\phi} = \argmin_{u \in \mathbb{R}^n} \textnormal{E}_{x|y}[D_\phi(u,x)]$ as follows
\begin{align*}
\hat{x}_{D_\phi} &= \argmin_{u \in \mathbb{R}^n} \phi(u) + \textnormal{E}_{x|y}[\phi(x)] - u^\top \textnormal{E}_{x|y}[\nabla \phi(x)] - x^\top \textnormal{E}_{x|y}[\nabla \phi(x)],\\
&=\argmin_{u \in \mathbb{R}^n} \phi(u) - u^\top \textnormal{E}_{x|y}[\nabla \phi(x)].
\end{align*}
In a manner akin to \cite{Lucka:2014}, the proof is concluded by using the divergence theorem, together with the fact that $p(x|y)$ is continuous and vanishes at least exponentially as $\|x\| \rightarrow 0$, to show that $\textnormal{E}_{x|y}[\nabla\phi(x)] = \int_{\mathbb{R}^n} \nabla p(x|y)\textnormal{d}x = 0$. Hence,
\begin{align*}
\hat{x}_{D_\phi} &= \argmin_{u \in \mathbb{R}^n} \phi(u),\\ 
&= \hat{x}_{\text{MAP}}. 
\end{align*}
Note that in the case where $p(x|y)$ involves hard constraints on the parameter space then generally $\textnormal{E}_{x|y}[\nabla\phi(x)] \neq 0$, and we have $\hat{x}_{D_\phi} = \argmin_{u \in \mathbb{R}^n} D_\phi(u,\hat{x}_{\text{MAP}}) - u^\top \textnormal{E}_{x|y}[\nabla\phi(x)]$ generally different from $\hat{x}_{\text{MAP}}$.

Finally, the proof of the third part of Theorem \ref{Theo2} follows directly from \cite[Proposition 1]{Banerjee:2005}, which for completeness we detail below
\begin{align*}
\hat{x}_{D^*_\phi} &= \argmin_{u \in \mathbb{R}^n} \textnormal{E}_{x|y}[D^*_\phi(u,x)],\\ 
&= \argmin_{u \in \mathbb{R}^n} \textnormal{E}_{x|y}[D_\phi(x,u)],\\
&= \argmin_{u \in \mathbb{R}^n} \textnormal{E}_{x|y}[D_\phi(x,u)] - \textnormal{E}_{x|y}[D_\phi(x,\hat{x}_{\text{MMSE}})],\\
& = \argmin_{u \in \mathbb{R}^n} \phi(\hat{x}_{\text{MMSE}}) - \phi(u) - (\hat{x}_{\text{MMSE}}-u)^\top \nabla\phi(u),\\
&= \argmin_{u \in \mathbb{R}^n} D_\phi(\hat{x}_{\text{MMSE}},u),\\
&=\hat{x}_{\text{MMSE}}.
\end{align*}

\noindent{Notice that strict log-concavity suffices to prove the second and third parts of Theorem \ref{Theo2}.}

\subsection*{Proof of Proposition \ref{Theo3}}
Assume that $\phi(x) = -\log p(x|y)$ is convex on $\mathbb{R}^n$ and $\mathcal{C}^1$. From the optimality condition of $\hat{x}_{\text{MAP}}$, $\nabla\phi(\hat{x}_{\text{MAP}}) = 0$ and hence the dual divergence
$$
D^{*}_{\phi}(\hat{x}_{\text{MAP}},x) = \phi(x) - \phi(\hat{x}_{\text{MAP}}) \, .
$$
Noting that $\textnormal{E}_{x|y}\left[{\phi(x)}\right] $ is the entropy of $x|y$, we use Proposition I.2 of \cite{Bobkov2011b} and obtain
$$
\textnormal{E}_{x|y}\left[{D^{*}_{\phi}(\hat{x}_{\text{MAP}},x)}\right]  = \textnormal{E}_{x|y}\left[{\phi(x)}\right] - {\phi(\hat{x}_{\text{MAP}})}\leq n.
$$
Finally, {using that $\hat{x}_{\text{MMSE}}$ minimises the posterior expectation of $D^{*}_{\phi}(\hat{x}_{\text{MMSE}},x)$ we obtain}
$$
\textnormal{E}_{x|y}\left[{D^{*}_{\phi}(\hat{x}_{\text{MMSE}},x)}\right] \leq \textnormal{E}_{x|y}\left[{D^{*}_{\phi}(\hat{x}_{\text{MAP}},x)}\right] \leq n\, ,
$$
concluding the proof.

\end{document}